\documentclass[10pt]{amsart}
\usepackage{latexsym}
\usepackage{amssymb}
\usepackage{amsmath}
\usepackage{mathrsfs}
\usepackage{calrsfs}
\usepackage{bbold}

 \def\Orth{\mathop{\fam0 Orth}}

 \def\upwardarrow{\mathord{
  \hbox to 5pt{\hss$\vcenter{\hbox to 2.4pt{\hss$\mathchar"222$\hss}\hrule}\hss$}
}}
\def\downwardarrow{\mathord{
  \hbox to 5pt{\hss$\vcenter{\hrule\hbox to 2.4pt{\hss$\mathchar"223$\hss}}\hss$}
}}

 \def\mix{\mathop{\fam0 mix}\nolimits}

 \def\[{\mathopen{\kern1pt\/ \vrule height7.6pt depth2.2pt width1pt\kern1.5pt}}
\def\]{\mathclose{\kern1.5pt\vrule height7.6pt depth2.2pt width1pt\kern1pt}}
\def\Bnorml{\mathopen{\kern1pt\/ \vrule height17pt depth10.2pt width1pt\kern1.5pt}}
\def\Bnormr{\mathclose{\kern1.5pt \vrule height17pt depth10.2pt width1pt\kern1pt}}
\def\bnorml{\mathopen{\kern1pt\/ \vrule height13pt depth7.2pt width1pt\kern1.5pt}}
\def\bnormr{\mathclose{\kern1.5pt \vrule height13pt depth7.2pt width1pt\kern1pt}}
\newcommand{\Lemma}[1]{\smallskip\textbf{{Lemma#1}.\/~}\sl}

\newcommand{\Theorem}[1]{\smallskip\textbf{{Theorem#1}.\/ }\sl}

\newcommand{\Corollary}[1]{\smallskip\textbf{{Corollary#1}.\/~}\sl}

\newcommand{\Definition}[1]{\smallskip{\sc Definition#1.}}
\newcommand{\Remark}[1]{\smallskip{\sc Remark#1.}~}

\newcommand{\Endproc}{\rm}

 \def\Orth{\mathop{\fam0 Orth}}

 \def\fin{\mathop{\fam0 fin}}

 \def\tv@rt{{\vert\mkern-2mu\vert\mkern-2mu\vert}}
 \def\tvert#1\tvert{\mathopen\tv@rt#1\mathclose\tv@rt}

\begin{document}

 \title
[Stone--Weierstrass Approximation Revisited]
{
Stone--Weierstrass Approximation Revisited
}

 \author{A.~G.~Kusraev and S.~S.~Kutateladze}

 \begin{abstract}
 The aim of the present article is to extend the
 Stone--Weierstrass
  theorem to  functions ranging in a lattice normed space and order rather than topological approximation.
We proceed with the  machinery of Boolean valued transfer from lattice normed space to normed space.
\end{abstract}


\maketitle

\section{Introduction}
\label{Intro}

   There are many generalizations of the classical Stone--Weierstrass approximation
 theorem which broaden the class of the continuos  scalar or vector valued  functions
.to be approximated  One of the most striking results belongs to Bishop; see \cite{Bi}.
The available  proof of Bishop's result utilizes  some nontrivial tools of functional analysis.
 Machado  formulated a version of the theorem for vector-valued functions and
 gave a completely elementary proof in \cite{Mac}. Ransford \cite{Ran} succeeded in finding the proof
that is both elementary and concise. Prolla  revised Machado's proof in~\cite{Pro} and found some new forms
 of the Stone-Weierstrass theorem that are given below.

 Throughout this paper, $S$ denotes a compact Hausdorff topological space,
\textit{compact
 space\/} for short. Let $X$ be a real or complex normed space and let $C(S,X)$
br the vector space
 of all continuous functions from $S$ to $X$ with the supremum norm
 $$
 \|f\|\!:=\sup\{\|f(s)\|:\ s\in S\}\quad(f\in C(S,X)).
 $$
 A~\textit{multiplier} of a subset $W\subset C(S,X)$ is a continuous function 
$\varphi:S\to[0,1]$ such that $\varphi f+(1_S-\varphi)g\in W$ for all $f,g\in W$, 
where $1_S$ is the identically $1$ function on $S$. The next 
three results correspond respectively to Theorems 1, 2, and 3 in~\cite{Pro}.

 \Theorem{~1.1}~Let $W$ be a nonempty subset of $C(S,X)$ such that the set of
all multipliers of~$W$ separates the points of\/ $S$. Then for ally $f\in C(S,X)$ and
 $0<\varepsilon\in\mathbb{R}$ the following  are equivalent:

 $(1)$~There exists $g\in W$ with $\|f-g\|<\varepsilon$.

 $(2)$~For every $s\in S$ there exists $g_s\in W$ such that $\|f(s)-g_s(s)\|_X<\varepsilon$.
 \Endproc

 \Theorem{~1.2}Let $W$ be a nonempty subset of $C(S,X)$ such that the set of all
 multipliers of\/ $W$ separates the points of~$S$. Then for every $f\in C(S,X)$
 there exists $s\in S$ such that
 $$
 \inf\{\|f-g\|:\ g\in W\}=\inf\{\|f(s)-x\|_X:\ x\in W(s)\}. $$
 \Endproc

 \Theorem{~1.3}Let $W$ be a vector subspace of $C(S,X)$ and
 $A\!:=\{\varphi\in C(S,\mathbb{R}):\ (\forall\,g\in W)\,\varphi g\in W\}$. Assume that $A$
 separates the points of $S$ and $W(s)$ approximates every member of $X$ for all $s\in S$.
 Then~$W$ is uniformly dense in $C(S,X)$.
 \Endproc

 The aim of the present article is to extend the above results to
the order rather than topological approximation of functions
ranging in a~lattice normed space.

 We use the standard notation and terminology of Aliprantis and Burkinshaw \cite{AB} for
 the theory of vector lattices. Everywhere below $E$ and $F$ denote some  Archimedean real
vector lattices, while $\mathbb{B}(E)$ and $\mathbb{P}(E)$ stand for the Boolean algebras
of all bands  and band projections in $E$. Recall also
 that $x,y\in F$ are  \textit{disjoint\/} (in symbols $x\perp y$) if
 $|x|\wedge|y|=0$. Also, $M^\perp\!:=\{x\in F:\ (\forall\,y\in\in M)\,x\perp y\}$.
 Given a~vector lattice $E$, we denote, by $E^\delta$ and $E^u$ the Dedekind
 completion and the universal completion of a vector lattice $E$; moreover, we assume
that $E\subset E^\delta\subset E^u$.

  Throughout the sequel $\mathbb{B}$ is a complete Boolean algebra with join $\vee$, meet  
$\wedge$, complement $(\cdot)^\perp$, unit (top) $\mathbb{1}$ and zero (bottom) $\mathbb{0}$.
  A~\textit{partition of an element} $b\in \mathbb{B}$ in a Boolean algebra $\mathbb{B}$ is
  a family $(b_\xi)$ in $\mathbb{B}$ such that $b_\xi\wedge b_\eta=\mathbb{0}$ for all $\xi\ne\eta$
  and $b=\sup_\xi b_\xi$; if $b$ is the unity of $\mathbb{B}$ then $(b_\xi)$ is
a~\textit{partition of unity}. We let $:=$ denote the assignment by definition, while
  $\mathbb{N}$ and $\mathbb{R}$ symbolize the naturals and the reals.

\section{Preliminaries}
\label{Rrel}
 In what follows, we will need some information about lattice normed spaces, approximating
 sets, extreme operators, and lattice homomorphisms.
 \medskip

 \textsc{Lattice normed spaces.}~We use the abbreviation LNS for ``lattice normed space.''
 Consider a~(real or complex) vector space $X$ and a~real vector lattice $E$.

 \Definition{~2.1}~An~{\it $E$-valued norm\/} on $E$ is a~mapping $\[\cdot\]:X\to E_+$ such
 that~$\[x\]=0$ implies $x=0$, $\[\lambda x\]=|\lambda|\[x\]$, and $\[x+y\]\leq \[x\]+\[y\]$
 for all $x,y\in X$, and $\lambda\in\mathbb{R}(\mathbb{C})$.  The pair $(X,\[\cdot\])$ is
 called an~\textit{LNS over\/} $E$. An~$E$-valued norm $\[\cdot\]$ on $X$ as well as $X$ itself
 is $d$-{\it decomposable\/} if, for each decomposition $\[x\]=e_1+e_2$ with
 disjoint $e_1,e_2\in E_+$ and $x\in X$, there exist $x_1,x_2 \in X$ such that $x=x_1+x_2$
 and $\[x_k\]=e_k$ for all $(k\!:= 1,2)$.

 If $X$ is a $d$-decomposable LNS over $E$, then there is a natural mapping that associates
 to each $\pi\in\mathbb{P}(E_0)$, with   $E_0=\{\[x\]:\ x\in X\}^{\perp\perp}$, the linear
 projection $\hat{\pi}$ in $X$. The set $\mathbb{P}(X)\!:=\{\hat{\pi}:\ \pi\in\mathbb{P}(E)\}$
  with the order relation defined by letting $\pi\leq\rho$ if and only if
 $\pi\circ\rho=\pi$ is a Boolean algebra and the mapping $\pi\mapsto\hat{\pi}$ is a Boolean
 isomorphism from $\mathbb{P}(E)$ onto $\mathbb{P}(X)$. Moreover,
 \begin{equation}\label{D}
 \[\hat{\pi} x+\hat{\pi}^\perp y\]=\pi\[x\]+\pi^\perp\[y\]\quad(\pi\in\mathbb{P}(E);
 x,y\in X).
 \end{equation}
 In the sequel, we will assume that  $\{\[x\]:\ x\in X\}^{\perp\perp}=E$  for every LNS $X$ over $E$
and the two Boolean algebras $\mathbb{P}(E)$ and
 $\mathbb{P}(X)$ are identified; see (see \cite[\S~2.1]{DOP}).

 \Definition{~2.2}~We say that a~net $(x_\alpha)_{\alpha\in {\rm A}}$ in $X$
 {\it norm $o$-converges} or  {\it norm $r$-converges with regulator
 $e\in E_+$\/}) to $x\in X$ and write $x=no$-$\lim x_\alpha$  or, respectively,
 $x=nr$-$\lim x_\alpha$) if $o$-$\lim_\alpha\[x-x_\alpha\]=0$ or, respectively,
 $r$-$\lim_\alpha\[x-x_\alpha\]=0$ with regulator $e\in E_+$). A~net $(x_\alpha)$ is
 {\it norm  $o$-fundamental} or {\it norm $r$-fundamental with regulator
 $e\in E_+$} if the~net $(x_\alpha-x_\beta)_{(\alpha,\beta)\in{\rm A}\times{\rm A}}$
 norm $o$-vanishes  or norm $r$-vanishes with regulator $e\in E_+$. The
 set of the $no$-limits  or $nr$-limits of all $no$-convergent  or $nr$-convergent
 nets in $X$  comprised of the elements of some subset $X_0\subset X$ is  the
{\it $no$-closure} or {\it $nr$-closure} of $X_0$.

 \Definition{~2.3}~A~lattice-normed space $X$ is {\it norm  $o$-complete} or \textit{norm $r$-complete\/} if every
 norm $o$-fundamental  or, respectively, norm $r$-fundamental net in $X$ norm $o$-converges
 or norm $r$-converges to an~element of $X$. A~{\it Ba\-nach--Kantorovich
 space\/} over a vector lattice $E$ is a vector space $X$ with a decomposable norm
 $\[\cdot \]:X\rightarrow E$ which is \textit{norm $o$-complete}. A~Ba\-nach--Kantorovich
 space over $E$ is  \textit{universally complete\/} if the vector lattice
 $E$ is universally complete. Given two LNSs $X$ and $Y$ over the same vector lattice,
 an operator $T:X\to Y$ is  \textit{isometric\/} if $\[Tx\]=\[x\]$ for all $x\in X$.

 \Definition{~2.4}~By a \textit{universal completion} or a \textit{norm completion} of
an~LNS $(X,E)$ we mean a universally complete Banach--Kantorovich space $(X^u,E^u)$
or, respectively, a Banach--Kantorovich space $(\tilde{X},E^\delta))$ together with a linear
 isometry $\iota:X\to X^u\,(\tilde{X})$ such that each universally complete subspace of
 $(X^u,E^u)$ or, respectively, each decomposable norm complete subspace of $(\tilde{X},E^\delta))$
 including  $\iota(X)$ coincides with $X^u\,(\tilde{X})$. The $d$-\textit{decomposable hull\/} of $X$ is the least
 $d$-decomposable subspace $\bar{X}\subset X$ including $\iota(X)$.

 \Lemma{~2.5}~For each LNS $X$ over a vector lattice $E$ the following  hold:

 $(1)$~There is a universal completion of $X$  unique to within linear isometry.

 $(2)$~There is a norm $o$-completion of $X$  unique to within linear isometry.

 $(3)$~There is a $d$-decomposable hull of $X$  unique to within an isometry.
 \Endproc

 \textsc{Approximating sets.}~We will briefly be recall the concept of approximating subset in an LNS 
which was introduced by Gutman in \cite{Gut} and turnes out very useful in the general theory
 of LNSs as well as in the study of disjointness preserving operators; see also \cite{Gut1}. 
Let $X$ be a LNS over a Dedekind complete vector lattice $E$.

 \Definition{~2.6}~If $(b_\xi)_{\xi \in \Xi}$ is a partition of unity in $\mathbb{P}(X)$
 and $(x_\xi)_{\xi \in \Xi}$ is a family in $X$, then  $x\in X$ with
 $b_\xi x_\xi=b_\xi x$ for all $\xi \in\Xi$ is a {\it mixing of $(x_\xi)$ by}
  $(b_\xi)$  denoted by $\mix_{\xi \in \Xi}b_\xi x_\xi$.
 The set of all mixings of arbitrary or finite families in $U$ is  the \textit{cyclic hull}
or, respectively, \textit{finitely cyclic hull\/}) of $U$ which is denoted by $\mix(U)$
or, respectively,  by $\mix_{\fin}(U)$. Say that $X$ is $d$-complete if, for each partition of unity $(\pi_\xi)$
in $\mathbb{P}(E)$ and norm bounded family $(x_\xi)$ ($\[x_\xi\]\leq e$ for all $\xi$ and some $e\in E$) in $X$, 
there is $x\in X$ with $x=\mix \pi_\xi x_\xi$. It easy to verify that the (finitely) cyclic hull of a~set $U$ 
is the smallest (finitely) cyclic set that includes $U$. It follows from (\ref{D}) that for a set $U$ 
to be finitely cyclic, it  suffices that $U$ contains the sums $\pi x+\pi^\perp y$ for all $x,y\in V$ 
and $\pi\in \mathbb{P}r(X)$.%

 Lemma 2.5 can be supplemented as follows: Each LNS has a $d$-completion $\stackrel{=}{X}$ unique up to linear isometry.

 \Definition{~2.7}~Let $U$ be a subset of an LNS $X$. We say that $U$ \textit{$($orderly$)$ approximates}  $x\in X$
if $\inf_{u\in U}\[x-u\]=0$. We say that $U$ \textit{$($orderly$)$ approximates\/}  $W\subset U$ if $U$ approximates
every  element of $W$. A subset of $X$ is  \textit{$($order$)$ approximating\/} if it approximates $X$.
The set $X_0$ is  $no$-\textit{dense\/}  or  $nr$-\textit{dense} in $X$ if every member of $X$ is
the $no$-limit ($nr$-limit) of some net in $X_0$.

 \Lemma{~2.8}~Let $X$ be a $d$-decomposable LNS over a Dedekind complete vector lattice $E$. Given a subset
 $U$ and an element $u$ of $X$, the following  hold:

 $(1)$~$U$ approximates $u$ if and only if $u$ is the $no$-limit in $\bar{X}$ of some net of elements from $\mix_{\fin}(U)$.%

 $(2)$~$U$ approximates $u$ if and only if $u$ is the $nr$-limit in of elements from $\stackrel{=}{X}$ of
 some net of elements of $\mix(U)$.

 $(3)$~If, moreover, $E$ has an order unity $\mathbb{1}$ then $U$ approximates $u$ if and only if $u$
is the $nr$-limit in $\stackrel{=}{X}$  with  regulator $\mathbb{1}$ of some net of elements of $\mix(U)$.
 \Endproc

 \begin{proof}~See  \cite[Propositions 1.3, 1.6, and 1.8]{Gut}.~\end{proof}

 \Lemma{~2.9}The following properties of a subset $U$ of an LNS $X$ are equivalent:%

 $(1)$~$U$ is an approximating subset of $X$.

 $(2)$~$\mix_{\fin}(U)$ is an approximating subset of $\bar{X}$.

  $(3)$~$\mix(U)$ is an approximating subset of $\stackrel{=}{X}$.

 $(4)$~$\mix(U)$ is norm $r$-dense in $\stackrel{=}{X}$.

 $(5)$~$\mix_{\fin}(U)$ is norm $o$-dense in $\bar{X}$.
 \Endproc

 \begin{proof}~See  \cite[Propositions 1.4, 1.7, and 1.9]{Gut}.~\end{proof}

 \medskip

 \textsc{Extreme operators and lattice homomorphisms.}~The question of recovering convex
 sets of operators from extreme points was raised in the well-known paper of Bonsall,
 Lindenstrauss, and Phelps \cite{BLP}. The general results on the geometric structure of the
 support sets of sublinear operators, as well as the facts of use below, can be found
in~\cite[Chapter~2]{SBD}.

 A \textit{sublinear operator\/} $p:X\to E$ from a vector space $X$ to a vector lattice $E$ is 
a~subadditive and positive homogeneous mapping; i.e.,  $p(x+y)\leq p(x)+p(y)$ and
 $p(rx)=rp(x)$ for all $x,y\in X$  and positive $r\in\mathbb{R}$. The collection of all linear 
operators from $X$ into $E$ dominated by $p$ is  the {\it support set} or the 
{\it subdifferential at zero} of $p$ and denoted by $\partial p$. In symbols,
 $$
 \partial p:=\{T\in L(X,E): (\forall x\in X)\,Tx\leq p(x)\},
 $$
 where $L(X,E)$ is the space of all linear operators from $X$ into $E$. If $A$ is a~subring 
and sublattice of the orthomorphism ring $\Lambda\!:=\Orth (E)$ and $p$ is a~$A$-sublinear operator
 $(p(\pi x+\rho y)\leq\pi p(x)+\rho p(y)$ for all $x,y\in X$ and $\pi, \rho\in A_+)$, then
 the members of $\partial (p)$ are automatically module homomorphisms; see~\cite[Theorem 2.3.15]{SBD}.

 Denote by $l_\infty (S, E)$ the set of all (order) bounded mappings from $S$ into $E$;
 i.e., $f\in l_\infty (S, E)$ if and only if $f:S\to E$ and
 $\{f(s):s\in S\}$ is order bounded in $E$. It is easy to verify that $l_\infty (S, E)$ with the coordinate-wise algebraic operations and order is a~Dedekind complete vector lattice if so is $E$. Moreover, $l_\infty (S, E)$ is a faithful module
 over $\Lambda$ with  multiplication $(\lambda,f)\mapsto \lambda f\!:=\lambda\circ f$
 $(\lambda\in\Lambda,f\in l_\infty (S, E))$.

 The operator $\varepsilon_{S, E}$ acting from $l_\infty (S, E)$ into $E$ by the rule
 \begin{equation*}
 \varepsilon_{S}\!:=\varepsilon_{S, E}:f\mapsto\sup \{f(s):s\in S\}\quad (f\in l_\infty(S, E))
 \end{equation*}
 is  the {\it canonical sublinear operator} given by $S$ and $E$.

 \Lemma{~2.10}For $e\in E$ denote by $\bar{e}$ the constant function $s\mapsto e$
 $(s\in S)$. Then
 $$
 \partial(\varepsilon_{S, E})=\{T\in L^+(l_\infty(S, E),E):\
 (\forall\,e\in E)\,T\bar{e}=e\},
 $$
 where $L^+(G,E)$ is the cone of all positive operators from $G$ to $E$.
 \Endproc

 \begin{proof}~See \cite[2.1.4]{SBD}.~\end{proof}

 \Definition{~2.11}~Given $s\in S$, the $\delta$-{\it function\/}
 $\delta_s:l_{\infty} (S, E)\to E$ is defined as $\delta_s: f\mapsto f(s)$. The mixing 
$o$-$\sum_{s\in S}\pi_s\delta_s$ of a~family $(\delta_s)_{s\in S}$ by a partition of unity $(\pi_s)_{s\in S}$ in $\mathbb{P}(E)$ 
is  an~$E$-valued {\it pure state} on $S$. (Here $o$-$\sum$ means the pointwise order sum.)

 Clearly, all $\delta$-functions and all pure states are extreme points of $\partial(\varepsilon_{S})$. Moreover, every extreme point of $\partial(\varepsilon_S)$
can  be approximated by pure states.

 \Lemma{~2.12}The set of extreme points of the support set $\partial(\varepsilon_{S})$ coincides 
with the set of lattice homomorphisms from $l_{\infty} (S,E)$ into~$E$ which belong to $\partial\varepsilon_{S}$. 
Moreover, all members of $\partial(\varepsilon_{S})$ are $A$-linear.
 \Endproc

 \begin{proof}~See \cite[Theorem 3]{Kut3} or \cite[Theorem 2.2.9]{SBD}.~\end{proof}

 \Lemma{~2.13}Each extreme point of the convex set $\partial(\varepsilon_S)$ is a~pointwise
 $r$-limit of a~net of pure states.
 \Endproc

 \begin{proof}~This fact was established in~\cite[The Main
 Theorem]{KK}; see \cite[Proposition 2.4.8 and Theorem 2.4.11]{SBD} for details.~\end{proof}

 \Lemma{~2.14}~Let $E$ and $F$ be two vector lattices with $F$ Dedekind complete. If $G$ is
 a majorizing vector sublattice of $E$ and $T:G\to F$ is a lattice homomorphism, then $T$
 extends to all of $E$ as a lattice homomorphism.
 \Endproc

 \begin{proof}~This is the well known Lipecki--Luxemburg--Schep Theorem see
 \cite[Theorem 2.29]{AB}.~\end{proof}

\section{Main Results}
\label{sec:3}

From now on, we will assume that $E$ is a Dedekind complete vector lattice and $X$ is
an~LNS over $E$, while $\mathbb{B}=\mathbb{P}(E)$ and $\Lambda\!:=\Orth(E)$. Then $\Lambda$
 is a Dedekind complete $f$-algebra (under composition) with unit $I\!:=I_E$; see
 \cite[Theorems 2.45 and 2.59]{AB}. Before stating the
 results, we will introduce the main object of study, i.e., the space of uniformly norm
 continuous functions with values in $X$.
 Let $(S,\mathcal{U})$ be a uniform
 space with  uniformity $\mathcal{U}$.

 \Definition{~3.1}~Let $X$ be an LNS over a Dedekind complete vector lattice $E$. 
A~vector-function $f:S\to X$ is  {\it uniformly order continuous} if
 $$
 \inf\limits_{U\in\mathcal{U}}\sup\{\[f(u_1)-f(u_2)\]: \ (u_1,u_2)\in U\}=0.
 $$

 This amounts to saying that $f$ is norm bounded on $S$ (i.e., there exists $e\in E$ with $\[f(s)\]\leq e$ for all $s\in S$) and 
if $\[f\]\!:=\sup\{\[f(s)\]:s\in S\}$, then for every $\[f\]\leq e\in E$ and if $0<\varepsilon\in\mathbb{R}$ there exists a partition of unity $(\pi_\alpha)_{\alpha\in\mathcal{U}}$ in $\mathbb{P}(X)$ such that $\pi_\alpha\[f(s_1)-f(s_2)\]\leq\varepsilon e$
for all $\alpha\in\mathcal{U}$ and
 $(s_1,s_2)\in\alpha$. We denote by $C_{uo}(S,X)$ the vector space of all uniformly
 order continuous mappings from $S$ into $X$ endowed with the $E$-valued norm
 $f\mapsto\[f\]$. Obviously, $C_{uo}(S,\Lambda)\subset l_\infty (S,\Lambda)$.

 For a compact space $S$ there exists exactly one uniformity $\mathcal{U}$ on the  $S$ that induces the original topology of $S$;
 moreover this uniformity is totally bounded and complete \cite[Theorems 8.3.13 and 8.3.16]{E}. Clearly, $C_{uo}(S,X)=C(S,X)$ is
the space of continuous function with values in a normed space $X$ whenever $E=\mathbb{R}$.

 The vector space $C_{ou}(S,\Lambda)$ will be considered with
 pointwise multiplication: given $\varphi,\psi\in C_{uo}(S,\Lambda)$, we put $\varphi\psi:s\mapsto\varphi(s)\psi(s)$ $(s\in S)$.
 We also equip $C_{uo}(S,X)$ with  multiplication by elements of  $C_{uo}(S,\Lambda)$ by putting $\varphi f:s\mapsto\varphi(s)f(s)$ $(s\in S)$.

 \Lemma{~3.2}If $X$ is a decomposable LNS or a Banach--Kantorovich space over $E$ then so
 is $C_{uo}(S,X)$. Moreover, $C_{uo}(S,\Lambda)$ is an $f$-algebra and $C_{uo}(S,X)$ is
 a module over an $f$-algebra $C_{uo}(S,\Lambda)$ with $\[\varphi f\]\leq\[\varphi\]\[f\]$
 for all $\varphi\in C_{uo}(S,\Lambda)$ and $f\in C_{uo}(S,X)$.
 \Endproc

 \begin{proof}~This is immediate from Theorem 4.9 and Lemma 5.3 and can be derived without
 Boolean valued analysis on using \cite[Theorem 2.2.3]{DOP}.~\end{proof}

 \Lemma{~3.3}Let $\varepsilon_{\!S,\Lambda}^c$ be the restriction  to $C_{uo}(S,\Lambda)$
 of the sublinear operator $\varepsilon_{\!S,\Lambda}:l_\infty(S,\Lambda)\to\Lambda$.
 Every lattice homomorphism in the convex set $\partial(\varepsilon_{\!S,\Lambda}^c)$ is
 a~pointwise $r$-limit of a~net of pure states.
 \Endproc

 \begin{proof}~Note that $C_{uo}(S,\Lambda)$ is a majorizing sublattice of 
$l_\infty(S,\Lambda)$ as for every $e\in E$ the function $\bar{e}:s\mapsto e$ belongs to $C_{uo}(S,\Lambda)$. 
It follows that a lattice homomorphism $\sigma:C_{uo}(S,\Lambda)\to\Lambda$ 
extends to a lattice homomorphism $\hat{\sigma}:l_\infty(S,\Lambda)\to\Lambda$ by
 Lemma 2.14. If $\sigma\in\partial(\varepsilon_{\!S,\Lambda}^c)$ then $\sigma(\bar{e})=e$ for all $e\in E$ and hence 
$\sigma\in\partial(\varepsilon_{\!S,\Lambda})$ by Lemma 2.10. It remains to observe that $\hat{\sigma}$ is 
an extreme point of $\partial(\varepsilon_{\!S,\Lambda})$ in view of Lemma 2.12 and
 appeal to Lemma 2.13.~\end{proof}

 \Definition{~3.4}~Let $W$ be a nonempty subset of $C_{uo}(S, X)$. A function
 $\varphi\in C_{uo}(S,\Lambda)$ with $0\leq\varphi(s)\leq I$ for all $s\in S$ is
  a \textit{multiplier\/} of $W$ if $\varphi f+(I-\varphi)g\in W$ for every pair
 of elements $f,g\in W$. The set of all multipliers of $W$ is denoted by $\mu(W)$.

 It is clear that if $\varphi,\psi\in\mu(W)$, then $I-\varphi\in W$ and $\varphi\psi\in\mu(W)$.

 Denote by $\Sigma\!:=\Sigma(S,\Lambda)$ the set of all extreme points of the convex set $\partial(\varepsilon_{\!S,\Lambda}^c)$. Each extreme point of $\partial(\varepsilon_{\!S,\Lambda}^c)$ extends to an extreme point of $\partial(\varepsilon_{\!S,\Lambda})$ by Milman's theorem for support sets (see \cite[Theorem 2.210]{SBD}). Thus, taking Lemmas 2.10 into account, we see that $\Sigma$ comprises all lattice homomorphisms $\sigma:C_{uo}(S,\Lambda)\to \Lambda$ with $\sigma(\bar{\lambda})=\lambda$, where $\bar{\lambda}:s\mapsto \lambda$ $(s\in S,\,\lambda\in\Lambda)$.
 Given $s\in S$, the
 restriction of the corresponding $\delta$-function to $C_{uo}(S,E)$ will be denoted by the same symbol $\delta_s$.  Clearly,
the mapping $s\mapsto\delta_{s}$ is an injection of $S$ into $\Sigma$. This is how  we will identify $S$ with the corresponding subset
of $\Sigma$. Given $u,v\in E$ and
 $\sigma,\tau\in\Sigma$, define
 \begin{equation}\label{V}
 [u=v]\!:=\bigvee\{b\in\mathbb{B}:\ bu=bv\},\quad [u\ne v]\!:=I_E-[u=v];
 \end{equation}
 \begin{equation}\label{VV}
 [\sigma=\tau]\!:=\bigvee\{b\in\mathbb{B}:\ b\sigma=b\tau\},
 \quad[\sigma\ne \tau]\!:=I_E-[\sigma=\tau].
 \end{equation}

 Let $\tilde{X}$ be the norm completion of an LNS $X$. The Boolean isomorphism
 $\pi\mapsto\hat{\pi}$ from $\mathbb{P}(E)$ onto $\mathbb{P}(X)$ can be extended to a monomorphism of the ring $\Lambda$ into the ring of endomorphisms of the additive group of $\tilde{X}$. Hence, $\tilde{X}$ admits a faithful module structure over $\Lambda$. Show that every uniformly order continuous function
 $f\in C_{uo}(S,X)$ extends  canonically to some function
 $\tilde{f}:\Sigma\to\tilde{X}$.

 \Definition{~3.5}~Given $f\in C_{uo}(S,X)$ and a pure state
 $\sigma=\sum_{s\in S}\pi_s\delta_s$ with a partition of unity $(\pi_s)_{s\in S}$ in
 $\mathbb{B}$, we put $\tilde{f}(\sigma)\!:=\sum_{s\in S}\pi_sf(s)$ and note that
 $\[\tilde{f}(\sigma)\]=\sum_{s\in S}\pi_s\[f(s)\]=\sigma(\varphi)\leq\[f\]$. By
 Lemma 2.13 an arbitrary member $\sigma$ of $\Sigma$ is the pointwise $r$-limit of
 some~net $(\sigma_\alpha)$ of pure states, and so we put
 $\bar{f}(\sigma)\!:=\lim_{\alpha}\bar{f}(\sigma_\alpha)$.

 The existence of limits in $\tilde{X}$ and the soundness of the above definition follow from Lemma 5.3.
 Moreover, $\[f\]=\[\tilde{f}\]\!:=\sup\{\[\tilde{f}(\sigma)\]_{\tilde{X}}:\
 \sigma\in\Sigma\}$ and $[\sigma=\tau]\leq[\tilde{f}(\sigma)=\tilde{f}(\tau)]$ for
 all $\sigma,\tau\in\Sigma$, as can be easily seen from Lemmas 5.3 and 5.4. 
In the sequel, we will write $f(\sigma)$ instead of $\tilde{f}(\sigma)$.

 \Lemma{~3.6}~Given $\sigma,\tau\in\Sigma$ and $\mathcal{F}\subset C_{uo}(S,\Lambda)$, the
 following are equivalent:%

 $(1)$~$[\sigma\ne\tau]\leq\bigvee\{[\sigma(\varphi)\ne\tau(\varphi)]:\
 \varphi\in\mathcal{F}\}$.

 $(2)$~$[\sigma=\tau]\geq\bigwedge\{[\sigma(\varphi)=\tau(\varphi)]:\
 \varphi\in\mathcal{F}\}$.

 $(3)$~If $b\sigma\ne b\tau$ for some $b\in\mathbb{B}$ then there exists $\varphi\in\mathcal{F}$ 
such that $b\sigma(f)\ne b\tau(f)$.%

 $(4)$~For every $b\in\mathbb{B}$ we have $b\sigma=b\tau$ whenever
 $b\sigma(\varphi)=b\tau(\varphi)$ for all $\varphi\in\mathcal{F}$.
 \Endproc

  \begin{proof}~The equivalence $(1)\longleftrightarrow(2)$ is an immediate consequence of
  the infinite De Morgan laws in $\mathbb{B}$, while $(3)\longleftrightarrow(4)$ follows easily
  from the logical equivalence of $A\longrightarrow B$ and $\neg B\longrightarrow\neg A$. If there is $b\in\mathbb{B}$, $b\sigma\ne  b\tau$, and $(1)$ holds; then $0<b\leq[\sigma\ne\tau]$
  and hence $b_0\!:=b\wedge[\sigma(\varphi)\ne\tau(\varphi)]>0$ for some $\varphi\in\mathcal{F}$.
  It follows that $b_0\sigma(f)\ne b_0\tau(f)$ is true, but then so is $b\sigma(f)\ne b\tau(f)$
  which implies $(1)\longrightarrow(3)$. If (4) holds and $b\sigma(\varphi)=b\tau(\varphi)$
  or, equivalently, $b\leq[\sigma(\varphi)=\tau(\varphi)]$ for all $\varphi\in\mathcal{F}$,
  then $[\sigma=\tau]$, whence $(4)\longrightarrow(2)$.~\end{proof}

 \Definition{~3.7}~Say that a subset $\mathcal{F}\subset C_{uo}(S,\Lambda)$ \textit{separates the points of\/} $\Sigma$ if, given any points $\sigma,\tau\in\Sigma$,
 one and, hence, all conditions $(1)$--$(4)$ of Lemma 3.6 are satisfied.

 Now we have all prerequisites to formulating the main results.

 \Theorem{~3.8}~Let $X$ be an LNS over a Dedekind coplete vector lattice, let $\tilde{X}$ be
the  norm completion if $X$, and let $W$ be a nonempty subset of $C_{uo}(S,X)$. If the set $\mu(W)$ of
 all multipliers of\/ $W$ separates the points of $\Sigma$; then, for every $f\in C_{uo}(S,X)$, 
the following are equivalent:%

 $(1)$~$W$ approximates $f$ in $C_{uo}(S,X)$.

 $(2)$~$W(\sigma)\!:=\{g(\sigma)\in\tilde{X}:\ g\in W\}$ approximates $f(\sigma)$ in $\tilde{X}$
 for all $\sigma\in\Sigma$.
 \Endproc

 \Corollary{~3.9}Let $X$ be a Banach--Kantorovich space and let $W$ be a nonempty subset of $C_{uo}(S,X)$ and $f\in C_{uo}(S,X)$. If $\mu(W)$ separates the points of $\Sigma$, then the following are equivalent:

 $(1)$~$f$ is the norm $o$-limit in $C_{uo}(S,X)$ of some net in $\mix_{\fin}(W)$.

 $(1')$~$f$ is the norm $r$-limit in $C_{uo}(S,X$ of some net from $\mix(W)$.

 $(2)$~For every $\sigma\in\Sigma$ the value $f(\sigma)$ is the norm $o$-limit in $ {X}$ of some net 
in~$\mix_{\fin}(W(\sigma))$.

 $(2')$~For every $\sigma\in\Sigma$ the value $f(\sigma)$ is the norm $r$-limit in $ {X}$ of some
 net from t $\mix(W(\sigma))$.
 \Endproc

 \begin{proof}~This is immediate from Theorem 3.8 and Lemma 2.8.~\end{proof}

 \Remark{~3.10}~It is worth highlighting the two extreme cases of Theorem 3.8 and Corollary 3.9:

 $(1)$~If $E=\mathbb{R}$, then $\mathbb{P}(X)=\{0,I_X\}$, $X=\bar{X}=\stackrel{=}{X}$, $\mix(U)=U$ for every
 $U\subset X$, and $S=\Sigma$. In this case, there is no need to involve $\tilde{X}$, since the extension operator $f\mapsto \tilde{f}$ in Definition 3.5 is the only place where the completeness of $\tilde{X}$ is needed. Thus, we may assume that $X$ is a normed space, and hence we arrive at Prolla's result; see \cite[Theorem 1 and Corollary 1]{Pro}.

 $(2)$~Another extreme case is $X=E$ and $\[x\]=|x|$ for all $x\in E$. In this event $\tilde{f}(\sigma)=\sigma(f)$ and we obtain a new version of the Stone--Weierstrass theorem for vector-functions with values in a Dedekind complete vector lattice:
 \textsl{If\/ $W$ is a nonempty subset of $C_{uo}(S,E)$ such that  $\mu(W)$ separates the points of\/ $\Sigma$, then for every $f\in C_{uo}(S,E)$ the following are equivalent:}

  $(i)$~$\inf\{\[f-g\]:\ g\in W\}=0$ \textsl{in\/} $C_{uo}(S,E)$.

 %

 $(ii)$~$\inf\{|\sigma(f)-\sigma(g)|:\ g\in W\}=0$ \textsl{in $E$ for all $\sigma\in\Sigma$}.

 Let us formulate this result with approximation  in terms of order convergence.

 \Corollary{~3.11}~Let $W$ be a nonempty subset of $C_{uo}(S,E)$ such that the set $\mu(W)$
 of all multipliers of\/ $W$ separates the points of\/ $\Sigma$. Then for each
 $f\in C_{uo}(S,E)$ the following are equivalent:

 $(1)$~$f$ is the norm $o$-limit in $C_{uo}(S,E)$ of some net in $\mix_{\fin}(W)$.

 $(2)$~For every $\sigma\in\Sigma$ the value $\sigma(f)$ is the $o$-limit in $E$ of some net in $\mix_{\fin}(\sigma(W))$, where $\sigma(W)=\{\sigma(g)\in E:\ g\in W\}$.
 \Endproc

 \Theorem{~3.12}~Let $W$ be a nonempty subset of $C_{uo}(S,X)$ such that the set $\mu(W)$
 of all multipliers of\/ $W$ separates the points of\/ $\Sigma$. Then for each
 $f\in C_{uo}(S,X)$ there exists $\sigma\in\Sigma$ such that
 $$
 \inf\{\[f-g\]:\ g\in W\}=\inf\{\[f(\sigma)-x\]_{\tilde{X}}:\ x\in W(\sigma)\}.
 $$
 \Endproc

  \Theorem{~3.13}Let $X$ be a Banach--Kantorovich space, let $W$ be a vector subspace of
 $C_{uo}(S,X)$ and $A\!:=\{\varphi\in C_{uo}(S,\Lambda):\ (\forall\,g\in W)\varphi g\in W\}$.
 Assume that $A$ separates the points of $\Sigma$ and $W(\sigma)$ approximates every member of $X$ for all $\sigma\in\Sigma$. Then the following hold:

 $(1)$~$\mix_{\fin}(W)$ is norm $o$-dense in $C_{uo}(S,X)$;

 $(3)$ if, moreover, $E$ has an order unity $\mathbb{1}$;, then $\mix(W)$ is norm
 $r$-dense in $C_{uo}(S,X)$ with the same regulator $\mathbb{1}$.
 \Endproc

  \begin{proof}~Demonstration may proceed along the lines of Theorems 3.8 and 3.12
 by Boolean valued analysis. However, for diversity, we will deduce Theorem 3.13. from Corollary
 3.9 by the standard means. Observe first that $A$ is a subalgebra of $C_{uo}(S,\Lambda)$
 and $\mu(W)=\{\varphi\in A:\ 0\leq\varphi\leq I\}$, where $I\!:=I_E\in\Lambda$. 
Prove that $\mu(W)$ separates the points of $\Sigma$. Given $\sigma,\tau\in\Sigma$ and $\varphi\in A$, put
 $\varphi_0(\omega)\!:=\varphi(\omega)-\varphi(\sigma)$ and
 $\psi(\omega)\!:=a\varphi_0^2(\omega)$ for all $\omega\in\Sigma$, where $a$ is the unique
 member of $\Lambda^u$ such that $a\[\varphi_0^2\]=\pi$ with $\pi=[a]=\big[\[\varphi_0^2\]\,\big]$ and $[a]$ being the order projection onto the band $\{a\}^{\perp\perp}$. Clearly,
$\psi\in A$, $\psi(\sigma)=0$ and
 $0\leq\psi(\omega)\leq\[a\varphi_0^2\]=a\sup\{|\omega(\varphi_0^2)|:\ \omega\in\Sigma\}
 =a\[\varphi_0^2\]=\pi\leq I_E$. It follows that $\psi\in\mu(W)$. Moreover,
 \begin{multline*}
 \ \ [\varphi(\tau)=\varphi(\sigma)]=[\varphi_0(\tau)=0]=[\varphi_0^2(\tau)=0]
 \\
 =[a\varphi_0^2(\tau)=0]=[\psi(\tau)=0]=[\psi(\tau)=\psi(\sigma)].\ \ \
 \end{multline*}
 Hence, for every $\varphi\in A$ there exists $\psi\in\mu(W)$ such that
 $[\varphi(\tau)=\varphi(\sigma)]=[\psi(\tau)=\psi(\sigma)]$. By hypothesis, $A$ separates the points of $\Sigma$, which means by Lemma 3.6 that $[\sigma=\tau]\geq\bigwedge\{[\sigma(\varphi)=\tau(\varphi)]:\ \varphi\in A\}$. Consequently, $[\sigma=\tau]\geq\bigwedge\{[\sigma(\psi)=\tau(\psi)]:\ \varphi\in\mu(W)\}$ and $\mu(W)$ separates the points of $\Sigma$. Denote by $\bar{W}$ the norm $o$-closure of $\mix_{\fin}(W)$ in $C_{uo}(S,X)$. By Corollary 3.9, $\bar{W}=C_{uo}(S,X)$. It remains to note that the norm $o$-closure of $\mix_{\fin}(W)$ coincide with the norm $r$-closure of $\mix(W)$ and, if there is an order unit $\mathbb{1}$ in $E$, the norm $r$-closure can be taken with respect to the same regulator $\mathbb{1}$.~\end{proof}

\section{Boolean Valued Requisites}
\label{sec:4}

 In the sequel, $\mathbb{B}$ is a complete Boolean algebra, $\mathbb{V}^{(\mathbb{B})}$
 is a corresponding \textit{Boolean valued model\/} of set theory. As the standard
 model of set theory, we consider the \textit{von Neumann universe\/} $\mathbb{V}$. We
 need some properties of $\mathbb{V}$ and $\mathbb{V}^{(\mathbb{B})}$ as well as some
 relationships between them; the detailed presentation can be found in \cite{Bell, BVA, KKTop}.

 There is a natural way of assigning to each statement $\phi$ about $x_1,\dots,x_n\in\mathbb{V}^{({\mathbb B})}$
the \textit{Boolean truth-value\/} $[\![\phi(x_1,\dots,x_n)]\!]\in\mathbb{B}$. The sentence $\phi(x_1,\dots,x_n)$ is called
 {\it true within $\mathbb{V}^{(\mathbb{B})}$} if $[\![\phi(x_1,\dots,x_n)]\!]=\mathbb{1}$.
 All axioms and rules of inference of the first-order predicate calculus with equality are true in $\mathbb{V}^{(\mathbb{B})}$.
 In particular,
 \begin{equation}\label{F}
 [\![u=v]\!]\wedge[\![\phi(u)]\!]\leq[\![(v)]\!]
 \end{equation}
 for all $u,v\in\mathbb{V}^{(\mathbb{B})}$ and every formula $\phi(x)$.
 It follows that all the theorems of $\rm ZFC$ (Zermelo--Fraenkel set theory with the axiom of choice) are true in $\mathbb{V}^{(\mathbb{B})}$.
 This statement is known as the Boolean valued \textit{transfer principle} or {\it transfer} for short.
There is also the \textit{maximum principle}, which enables us
to construct all  particular objects in $\mathbb{V}^{(\mathbb{B})}$.

 Moreover, there is a~smooth mathematical technique for interplay between the interpretations of any
  fact in the two models~${\mathbb{V}}$ and~${\mathbb{V}}^{({\mathbb B})}$.
 The relevant {\it ascending-and-descending machinery\/} rests on the functors of
 \textit{canonical embedding} (or \textit{standard name}\/)
 $X\mapsto X^{\scriptscriptstyle\wedge}$ and  \textit{ascent} $X\mapsto X{\uparrow}$,
 both acting from $\mathbb{V}$ to $\mathbb{V}^{({\mathbb B})}$, and the functor of
 \textit{descent} $X\mapsto X{\downarrow}$, acting from $\mathbb{V}^{({\mathbb B})}$
 to $\mathbb{V}$; see \cite{Bell, BVA, KKTop} for details.

Observe some simple properties of the standard name mapping we need in the
sequel:

 \Lemma{~4.1} 
 $(1)$ Given  $x\in \mathbb{V}$ and a~formula
 $\varphi$ of~$ZF$ $($Zermelo--Fraenkel set theory$)$, we have
 \begin{equation*}
 [\![(\exists\, y\in
 x^{\scriptscriptstyle\wedge})\,\varphi(y)]\!]=
 {\textstyle\bigvee}\{[\![\varphi(z^{\scriptscriptstyle\wedge})]\!]:\,
 z\in x \},
 \end{equation*}
 \begin{equation*}
 \ [\![(\forall\, y\in x^{\scriptscriptstyle\wedge})\,\varphi(y)]\!]
 ={\textstyle\bigwedge}\{[\![\varphi(z^{\scriptscriptstyle\wedge})]\!]:\,
 z\in x \}.
 \end{equation*}

 $(2)$~The standard name mapping is injective. Moreover,
 for all ${x,y\in\mathbb{V}}$ we have
\begin{equation*}
 x\in y\longleftrightarrow \mathbb{V}^{(\mathbb{B})}\models
 x^{\scriptscriptstyle\wedge}\in y^{\scriptscriptstyle\wedge},
 \end{equation*}
 \begin{equation*}
 x= y\longleftrightarrow \mathbb{V}^{(\mathbb{B})}\models
 x^{\scriptscriptstyle\wedge}=y^{\scriptscriptstyle\wedge}.
\end{equation*}
 \Endproc

 \Lemma{~4.2}Denote by $\mathrsfs P_{\fin}(X)$ the collection of all finite sunsets of $X\in\mathbb{V}$ and let
 $[\![\mathrsfs P_{\fin}(Y)$  the collection of all finite
 subsets of $Y]\!]=\mathbb{1}$ with $Y\in\mathbb{V}^{(\mathbb{B})}$. Then
 $$
 \mathbb{V}^{(\mathbb{B})}\models\mathrsfs P _{\fin} (X^{\scriptscriptstyle\wedge})=
 \mathrsfs P_{\fin} (X)^{\scriptscriptstyle\wedge}.
 $$
 \Endproc

 \Lemma{~4.3}Let $\varnothing\ne X\in\mathbb{V}$, $Y\in\mathbb{V}^{(\mathbb{B})}$, and
 $[\![Y\ne\varnothing]\!]=\mathbb{1}$. Denote by $\mathcal{F}(X,Y{\downarrow})$ the set of
 all functions from $X$ to $Y{\downarrow}$ and let
 $\mathcal{F}(X^{\scriptscriptstyle\wedge},Y){\downarrow}$ stand for the set of all
 functions from $X^{\scriptscriptstyle\wedge}$ to $Y$ within $\mathbb{V}^{(\mathbb{B})}$.
 Then the following hold:

 $(1)$~If $[\![g$ is a~function from~$X^{\scriptscriptstyle\wedge}$ to~$Y]\!]=\mathbb{1}$,
 then there exists a function $g\downwardarrow$ from $X$ to $Y{\downarrow}$ uniquely determined by
 $$
 [\![g\downwardarrow(x)=g(x^{\scriptscriptstyle\wedge})]\!]= \mathbb{1} \quad (x\in X).
 $$

 $(2)$~If $f$ is a~function from $X$ to~$Y{\downarrow}$, then there exists a function
 $f\upwardarrow$ from $X^{\scriptscriptstyle\wedge}$ to $Y$ within
 $\mathbb{V}^{(\mathbb{B})}$ determined uniquely by
 $$
 [\![f\upwardarrow(x^{\scriptscriptstyle\wedge})=f(x)]\!]=\mathbb{1} \quad(x\in X).
 $$

 $(3)$~The mappings $f\mapsto f\upwardarrow$ and $g\mapsto g\downwardarrow$ are  inverse
to one another  and establish  bijections between $\mathcal{F}(X,Y{\downarrow})$ and
 $\mathcal{F}(X^{\scriptscriptstyle\wedge},Y){\downarrow}$.

 $(4)$~$[\![f\upwardarrow(A^{\scriptscriptstyle\wedge})=f(A){\uparrow}]\!]=\mathbb{1}$
 for every $A\subset X$.
 \Endproc

 The above functors are applicable, in particular, to algebraic structures. Applying the ransfer and maximum principles to the $\rm ZFC$-theorem on the existence of the fields of reals, we will find $\mathcal{R}\in \mathbb{V}^{(\mathbb{B})}$, called the \textit{reals within\/} $\mathbb{V}^{(\mathbb{B})}$ satisfying $[\![\mathcal{R}$ is the reals$]\!]=\mathbb{1}$ and $[\![1^{\scriptscriptstyle\wedge}\in\mathbb{R}
 ^{\scriptscriptstyle\wedge}\subset\mathcal{R}]\!]=\mathbb{1}$, where
 $\mathbb{R}\in\mathbb{V}$ is the (standard) field of reals with unit $1$.
 \textit{Gordon's theorem\/} \cite{Gor} establishes the relationship between
 $\mathbb{R}$, $\mathcal{R}$, and $\mathcal{R}{\downarrow}$.

  \Theorem{~4.4 {\rm (Gordon, 1977)}}~The descent $\mathcal{R}{\downarrow}$ of $\mathcal{R}$
 $($with the descended operations and order$)$ is a universally complete vector lattice
 with weak order unit $\mathbb{1}\!:=1^{\scriptscriptstyle\wedge}$. Moreover, the field
 $\mathcal{R}\in\mathbb{V}^{(\mathbb{B})}$ can be chosen so that
 $[\![\mathbb{R}^{\scriptscriptstyle\wedge}$ is a dense subfield of $\mathcal{R}]\!]=\mathbb{1}$.
 \Endproc

 The detailed presentation of the proofs of Gordon's theorem and the following two
 corollaries can be found in \cite[Theorems 2.2.4 and 2.3.2]{KKTop}.

 \Corollary{~4.5}There is a Boolean isomorphism $\chi$ from $\mathbb{B}$ onto $\mathbb{P}(\mathcal{R}{\downarrow})$ such that for all $x,y\in\mathcal{R}{\downarrow}$ and $b\in\mathbb{B}$ we have
  \begin{equation}\label{BV}
 \chi (b) x=\chi (b) y\longleftrightarrow b\leq[\![\,x=y\,]\!],
 \end{equation}
 \begin{equation*}
 \chi (b) x\le \chi (b) y\longleftrightarrow b\le [\![\,x\leq y\,]\!].
 \end{equation*}
 \Endproc

 \Corollary{~4.6}The universally complete vector lattice~$\mathcal{R}{\downarrow}$
 with the descended multiplication is a~semiprime $f$-algebra with the order and ring
 unit~$\mathbb{1}\!:=1^{\scriptscriptstyle\wedge}$. Moreover, for every~$b\in\mathbb{B}$
 the band projection $\chi(b)$ acts as multiplication by the $\chi (b)\mathbb{1}$.
 \Endproc

 \Lemma{~4.7}The following equivalences hold for a~nonempty
 set~$A\subset \mathrsfs R{\downarrow}$ and all~$a\in \mathrsfs R$ and
 $b\in\mathbb{B}$:
 \begin{equation}\label{S} 
 b\le [\![\,a=\sup
 (A{\uparrow})\,]\!]\longleftrightarrow\chi (b)
 a =\sup \chi (b) (A),
 \end{equation}
 \begin{equation}\label{I}
 b\le [\![\,a=\inf
 (A{\uparrow})\,]\!]\longleftrightarrow\chi (b)
 a=\inf \chi (b) (A).
 \end{equation}
 \Endproc

 \Definition{~4.8}~Take a~normed space $\mathcal{X}\!:=(\mathcal X,\rho)$ within $\mathbb{V}^{(\mathbb{B})}$,
that is $[\![\rho:\mathcal{X}\to\mathcal{R}$ is a norm on a (real or complex) vector space $\mathcal{X}]\!]=\mathbb{1}$.
The {\it descent} $\mathcal{X}{\downarrow}$ of $\mathcal{X}$ is a pair $(\mathcal{X}{\downarrow},\[\cdot\])$,
where $\[{\cdot}\]\!:=\rho{\downarrow}({\cdot}): \mathcal{X}{\downarrow}\to\mathcal{R}{\downarrow}$ is the descent of the intermal norm $\rho$.

 If $\mathcal{X}$ is a~Banach space within $\mathbb{V}^{(\mathbb{B})}$, then the descent
 $X\!:=\mathcal{X}{\downarrow}$ is a universally complete Banach--Kantorovich space over
$\mathcal{R}{\downarrow}$; see \cite[Theorem 5.4.1]{DOP}.

\Theorem{~4.9}For every LNS $X$ over $E$ with $E=\[X\]^{\perp\perp}$
 and\/ $\mathbb{B}\!:=\mathbb{P}(E)$ there exists
a Banach space  $\mathcal{X}$ within $\mathbb{V}^{(\mathbb{B})}$  unnique to within a linear isometry
and called  the $($\textit{Boolean valued representation}$)$ of~$X$ whose descent $\mathcal{X}{\downarrow}$ is
the universal completion of $X$.
 \Endproc

 \begin{proof}
 The proof can be found in \cite[Theorem 8.3.2]{DOP} and \cite[Theorem 5.4.2]{BVA}.
 \end{proof}

\section{Proofs}
\label{sec:5}

 In this section $\mathbb{B}\!:=\mathbb{P}(E)$, and $\mathbb{V}^{(\mathbb{B})}$ is the corresponding Boolean valued model of set theory.

\Lemma{~5.1}For every compact space $S$ there exists a~compact space $\tilde{S}$ within $\mathbb{V}^{(\mathbb{B})}$ such that $\iota:S\to\tilde{S}{\downarrow}$
unique $($up to homeomorphism$)$ and such that
 $[\![\iota(S){\uparrow}$ is dense in $\tilde{S}]\!]= \mathbb{1}$, where the embedding
 $\iota:S\to\tilde{S}{\downarrow}$ is defined as
 $\iota:s\mapsto s^{\scriptscriptstyle\wedge}$ $(s\in S)$.
 \Endproc

 \begin{proof}~For a compact space $S$ there is exactly one uniformity $\mathcal{U}$ on~$S$ that 
induces the original topology of $S$; moreover, $\mathcal{U}$ is totally
 bounded and complete \cite[Theorems 8.3.13 and 8.3.16]{E}. Working within $\mathbb{V}^{(\mathbb{B})}$, we observe 
 that $\mathcal{U}^{\scriptscriptstyle\wedge}$ may
 fail to be a uniformity on $S^{\scriptscriptstyle\wedge}$. However,
 $S^{\scriptscriptstyle\wedge}$ becomes a uniform space when endowes with  the uniformity base
 $\mathcal{U}^{\scriptscriptstyle\wedge}$; this uniformity we will denote by the same
 symbol $\mathcal{U}^{\scriptscriptstyle\wedge}$. Define the embedding
 $\iota:S\to S^{\scriptscriptstyle\wedge}$ as $s\mapsto s^{\scriptscriptstyle\wedge}$
 $(s\in S)$. By transfer  every uniform space has exactly one
 (up to a uniform isomorphism) completion, \cite[Theorems 8.3.12]{E}. Denote by
 $(\tilde{S},\tilde{\mathcal{U}})$ the completion of this uniform space
 $(S^{\scriptscriptstyle\wedge},\mathcal{U}^{\scriptscriptstyle\wedge})$ within $\mathbb{V}^{(\mathbb{B})}$. 
This and the equation $[\![\iota(S){\uparrow}=S^{\scriptscriptstyle\wedge}]\!]=\mathbb{1}$ (see Lemma 4.3(4)) imply 
that $\iota(S){\uparrow}$~is dense in~$\tilde{S}$. Now, to ensure that $\tilde{S}$ is compact
 it suffices to prove that $(S^{\scriptscriptstyle\wedge},\mathcal{U}^{\scriptscriptstyle\wedge})$ is  totally bounded; see \cite[Corollary 8.3.17]{E}.
 The latter amounts to checking that for every
 $U\in \mathcal{U}^{\scriptscriptstyle\wedge}$ there is a finite subset $\theta\subset\mathcal{P}_{\fin}(S^{\scriptscriptstyle\wedge})$ that is $U$-dense in  $(S^{\scriptscriptstyle\wedge},\mathcal{U}^{\scriptscriptstyle\wedge})$. In symbols,
 $$
 [\![(\forall\,U\in\mathcal{U}^{\scriptscriptstyle\wedge})\,
 (\exists\,\theta\in\mathcal{P}_{\fin}(S^{\scriptscriptstyle\wedge}))\,
 S^{\scriptscriptstyle\wedge}\subset U(\theta)]\!]=\mathbb{1},
 $$
 where $\mathcal{P}_{\fin}(S^{\scriptscriptstyle\wedge})$ is the collection of all finite subsets of $S^{\scriptscriptstyle\wedge}$, $U(\theta)=\bigcup_{s\in\theta}U(s)$, and
 $U(s)\!:=\{t\in S:\ (s,t)\in U\}$. Appreciating the equalities
$\mathcal{P}_{\fin}(S^{\scriptscriptstyle\wedge})=\mathcal{P}_{\fin}(S)^{\scriptscriptstyle\wedge}$
 (see Lemma 4.2) and $U(\theta)^{\scriptscriptstyle\wedge}=
 U^{\scriptscriptstyle\wedge}(\theta^{\scriptscriptstyle\wedge})$ (see \cite[Theorem 3.1.5(2)]{BVA}),
the simple rules for calculating Boolean truth values for quantifiers over standard names
 (see Lemma 4.1(1)), we arrive at the equivalent statement
 $$
 \mathbb{1}=\bigvee\{[\![S^{\scriptscriptstyle\wedge}\subset
 U(\theta)^{\scriptscriptstyle\wedge}]\!]:\
 \theta\in\mathcal{P}_{\fin}(S^{\scriptscriptstyle\wedge})\}.
 $$
 Since the Boolean truth value
 $[\![S^{\scriptscriptstyle\wedge}\subset U(\theta)^{\scriptscriptstyle\wedge}]\!]$ can
 take only the two values $\mathbb{0}\in\mathbb{B}$ and $\mathbb{1}\in\mathbb{B}$ (Lemma 4.1(2)), we obtain
 another equivalent form of the desired statement: For every $U\in\mathcal{U}$ there
 is $\theta\in\mathcal{P}_{\fin}(S)$ with $S\subset U(\theta)$; i.e.,  $S$ is totally bounded. The last claim is true by our hypothesis.
 \end{proof}

 \Remark{~5.2}The internal compact space $\tilde{S}$ is often called the \textit{Boolean extension} of $S$. The general theory of Boolean extensions of uniform spaces was built by Gordon and Lyubetskii \cite{GL}. The Boolean valued interpretation of compactness gives rise to the notion of a \textit{cyclically compact space\/} in such a way that the descent
 $\tilde{S}{\downarrow}$ turns out to be cyclically compact; see \cite[\S~8.5]{DOP}. The equivalent concept of \textit{mix-compact\/} subset of an LNS space was treated in Gutman and Lisovskaya \cite{GuL}; see also \cite[2.12.B, 2.12.C, and 2.12.13]{KKTop}.

 Define $C(\tilde{S},\mathcal{X})\in\mathbb{V}^{(\mathbb{B})}$ to be the set of continuous
 functions $f:\tilde{S}\to\mathcal{X}$ within $\mathbb{V}^{(\mathbb{B})}$; i.e., $\big[\!\!\big[C(\tilde{S},\mathcal{X})$ is the space of continuous functions from
 $\tilde{S}$ to $\mathcal{X}$ with the supremum norm $\|f\|\!:=\sup_{s\in\tilde{S}}
 \|f(s)\|_{\mathcal{X}}$\,$\big]\!\!\big]=\mathbb{1}$.

 \Lemma{~5.3} Let $X$ be a Banach--Kantorovich space and let $\mathcal{X}$ be the Boolean valued
 representations of $X$. For every $f\in C_{uo}(S,X)$ there exists a unique
 $\tilde{f}\in\mathbb{V}^{(\mathbb{B})}$ such that
 $[\![\tilde{f}\in C(\tilde{S},\mathcal{X})]\!]= \mathbb{1}$ and
 $[\![\tilde{f}(u^{\scriptscriptstyle\wedge})=f(u)]\!]= \mathbb{1}$ for
 all $u\in S$. The mapping $f\mapsto\tilde{f}$ is an linear isometry $($in the sense of lattice
 valued norms$)$ from $C_{uo}(S,X)$ into $C(\tilde{S},\mathcal{X}){\downarrow}$.
 \Endproc

 \begin{proof}~Without loss of generality we may assume that $X=\mathcal{X}{\downarrow}$ for some Banach space 
within $\mathbb{V}^{(\mathbb{B})}$. Let $f$ be a uniformly order continuous function from $S$ to $X$. 
Let the function $f{\upwardarrow}:S^{\scriptscriptstyle\wedge}\to X$
 within $\mathbb{V}^{(\mathbb{B})}$ be defined as in Lemma 4.3(2). If $f{\upwardarrow}$ is uniformly continuous 
then $f{\upwardarrow}$ extends uniquely  to some continuous
 functions $\tilde{f}:\tilde{S}\to\mathcal{X}$. By  transfer  the restriction operator $g\mapsto g|_{S^{\scriptscriptstyle\wedge}}$ from $C(\tilde{S},\mathcal{X})$ to
 $C_{uo}(S^{\scriptscriptstyle\wedge},\mathcal{X})$ is an isometric lattice isomorphism of
 Banach lattices within $\mathbb{V}^{(\mathbb{B})}$. Then the Banach--Kantorovich lattices
 $C(\tilde{S},\mathcal{X}){\downarrow}$ and
 $C_{uo}(S^{\scriptscriptstyle\wedge},\mathcal{X}){\downarrow}$ are also isometrically
 lattice isomorphic. So it remains to prove that $C_{uo}(S,X)$ and
 $C_{uo}(S^{\scriptscriptstyle\wedge},\mathcal{X}){\downarrow}$ are isometrically lattice
 isomorphic. In virtue of Lemma 4.3(3) it  suffices to ensure that $f{\upwardarrow}$ is
 uniformly continuous if and only if $f$ is uniformly order continuous. Define
 $D_U,d_U\in\mathbb{V}$ and $\mathcal{D}_U, \delta_U\in\mathbb{V}^{(\mathbb{B})}$ by
 \begin{equation*}
 D_U\!:=\{\[f(s_1)-f(s_2)\]:\ (s_1,s_2)\in U\},\quad d_U\!:=\sup(D_U)\in E;
 \end{equation*}
 \begin{equation*}
  \mathcal{D}_U\!:=\{\|f{\upwardarrow}(s_1)-f{\upwardarrow}(s_2)\|:\
 (s_1,s_2)\in U^{\scriptscriptstyle\wedge}\},\quad
 \delta_U\!:=\sup\mathcal{D}\in\mathcal{R}.
 \end{equation*}
 Observe that $[\![D_U{\uparrow}=\mathcal{D}_U]\!]=\mathbb{1}$ by Lemma 4.3(4) and $[\![d_U{\uparrow}=\delta_U]\!]=\mathbb{1}$
by Lemma 4.7 (formula \eqref{S}). Applying Lemma 4.7 again (formula \eqref{I}), we conclude that
 $\inf(D_U)=0$ if and only if $[\![\inf(\mathcal{D}_U)=0]\!]=\mathbb{1}$, as required.~\end{proof}

 \Lemma{~5.4}Let $\tilde{\Sigma}\in\mathbb{V}^{(\mathbb{B})}$ be defined as $\tilde{\Sigma}\!=\{\delta_s:\ s\in\tilde{S}\}$, where $\delta_s:C(\tilde{S},\mathcal{R})\to\mathcal{R}$ is the point evaluation $\varphi\mapsto\varphi(s)$ within $\mathbb{V}^{(\mathbb{B})}$ determined by a given $s\in\tilde{S}$. Then the mapping $j:s\mapsto(\delta_s){\downarrow}$ is a bijection
 from $\tilde{\Sigma}{\downarrow}$ onto $\Sigma$. Moreover, if
 $j(\tilde{\sigma})=\sigma$ and $j(\tilde{\tau})=\tau$ for some
 $\tilde{\sigma},\tilde{\tau}\in\tilde{\Sigma}$, then
 $[\![\tilde{\sigma}=\tilde{\tau}]\!]=[\sigma=\tau]$.
 \Endproc

 \begin{proof}~The mapping $\varkappa:\varphi\mapsto\tilde{\varphi}$ is a linear isometry $($in
 the sense of lattice valued norms$)$ from $C_{uo}(S,\mathcal{R}{\downarrow})$ onto
 $C(\tilde{S},\mathcal{R}){\downarrow}$; see Lemma 5.3. Moreover, $\varkappa$  is a lattice isomorphism
 in view of the easily verifiable fact that $\varphi\in C_{uo}(S,\mathcal{R}{\downarrow})$
 is positive if and only if so is $\tilde{\varphi}$ within $\mathbb{V}^{(\mathbb{B})}$.
 In this way we may identify $C_{uo}(S,\mathcal{R}{\downarrow})$ and
 $C(\tilde{S},\mathcal{R}){\downarrow}$. If $\sigma$ is a $\mathcal{R}{\downarrow}$-linear
 positive operator from $C_{uo}(S,\mathcal{R}{\downarrow})$ to $\mathcal{R}{\downarrow}$
 and $\tilde{\sigma}\!:=\sigma\circ\varkappa$, then $\tilde{\sigma}{\uparrow}$ is a positive
 linear functional on $C(\tilde{S},\mathcal{R})$ within $\mathbb{V}^{(\mathbb{B})}$. By transfer, $\tilde{\sigma}{\uparrow}\in\tilde{\Sigma}$ if and only if $\tilde{\sigma}{\uparrow}$ is a lattice homomorphism and
 $\tilde{\sigma}{\uparrow}(1_{\tilde{S}})=1$. It can be easily checked that the latter means that $\sigma$ is a lattice homomorphism and $\sigma(\bar{e})=e$ for all
 $e\in\mathcal{R}{\downarrow}$. The equality $[\![\tilde{\sigma}=\tilde{\tau}]\!]=[\sigma=\tau]$ is immediate from \eqref{BV}.
 \end{proof}

 \Lemma{~5.5}Let $E=\mathcal{R}{\downarrow}$ and $X=\mathcal{X}{\downarrow}$ with $\mathcal{X}$ 
being a Banach space within $\in\mathbb{V}^{(\mathbb{B})}$. Consider  $W\subset C_{uo}(S,X)$ 
and define $\widetilde{W}\in\mathbb{V}^{(\mathbb{B})}$ by the formula 
$\widetilde{W}\!:=\{\tilde{f}:\ f\in W\}{\uparrow}$. Then the following  hold:

 $(1)$~$[\![\widetilde{W}\subset C_{uo}(\tilde{S},\mathcal{X})]\!]=\mathbb{1}$ and
 $[\![\widetilde{\mu(W)}\subset\mu({\widetilde{W}})]\!]=\mathbb{1}$;

 $(2)$~if $\mu(W)$ separates the points of $\Sigma$, then $[\![\mu({\widetilde{W}})$ separates 
the points of $\tilde{S}]\!]=\mathbb{1}$.%
 \Endproc

 \begin{proof}
 $(1):$~The first inclusion is obvious. To prove the second, start with  observing
 that the mapping $\psi\mapsto\tilde{\psi}$ is a bijection from
 $C(\tilde{S},[0,1]){\downarrow}$ onto $C_{uo}(S,[0,\mathbb{1}])$. Indeed, if
 $\varphi=\tilde{\psi}$, then $\varphi(S^{\scriptscriptstyle\wedge})=\psi(S){\uparrow}$ by Lemma 4.3(4);
and, as $\varphi(\tilde{S})$ is the closure of $\varphi(S^{\scriptscriptstyle\wedge})$, we have $\varphi(\tilde{S})\subset[0,1]$ if and
 only if $\psi(S){\uparrow\downarrow}=\mix(\psi(S))\subset[0,\mathbb{1}]$ if and only if $\psi(S)\subset[0,\mathbb{1}]$.
 Now, using (\ref{F}) and the formulas for computing Boolean truth values for quantifiers over ascents (see \cite[1.6.2]{KKTop}), for every 
 $\varphi\in C(\tilde{S},[0,1]){\downarrow}$ we have
 \begin{align*}
 \Big[\!\!\Big[\varphi\in\widetilde{\mu(W)}\,\Big]\!\!\Big]
 =&\bigvee_{\psi\in\mu(W)}[\![\varphi=\tilde{\psi}]\!]
 \\
 =&\bigvee_{\psi\in\mu(W)}[\![\varphi=\tilde{\psi}]\!]\wedge
 [\![\tilde{\psi}\widetilde{W}+(1-\widetilde{\psi})\widetilde{W}\subset\widetilde{W}]\!]
 \\
 \leq& \bigvee_{\psi\in\mu(W)}[\![\varphi\widetilde{W}+(1-\varphi)\widetilde{W}\subset\widetilde{W}]\!]
 \\
 =&[\![\varphi\widetilde{W}+(1-\varphi)\widetilde{W}\subset\widetilde{W}]\!]
 \\
 =&[\![\varphi\in\mu(\widetilde{W})]\!].
 \end{align*}


 $(2):$~Assume that $\mu(W)$ separates the points of $\Sigma$. According to Lemma 3.6(1)
 this is equivalent to the equation $[\sigma\ne\tau]\Rightarrow\bigvee\{[\sigma(\varphi)
 \ne\tau(\varphi)]:\ \varphi\in\mu(W\}=\mathbb{1}$ for all $\sigma\tau\in\Sigma$. Note also
 that if $\varphi=\tilde{\psi}$, then $[\![\sigma{\uparrow}=\tau{\uparrow}]\!]=[\sigma=\tau]$
 and $[\![\sigma{\uparrow}(\varphi)=\tau{\uparrow}(\varphi)]\!]=[\sigma(\psi)=\tau(\psi)]$
 by (\ref{V}), (\ref{VV}), and (\ref{BV}).
 Denote $b\!:=[\![\widetilde{\mu(W)}$ separates
 the points of $\tilde{\Sigma}]\!]$. Easy calculations of Boolean truth values yield
 \begin{align*}
 b=&[\![(\forall\,\sigma,\tau\in\tilde{\Sigma})
 (\sigma\ne\tau)\rightarrow(\exists\varphi\in\widetilde{\mu(W)}
 \varphi(\sigma)\ne\varphi(\tau))]\!]
 \\
 =&\bigwedge_{\sigma,\tau\in\Sigma}\bigg([\![\sigma{\uparrow}\ne\tau{\uparrow}]\!]\Rightarrow
 \bigvee_{\psi\in\mu(W)}[\![\sigma{\uparrow}(\tilde{\psi})\ne\tau(\tilde{\psi})]\!]\bigg)
 \\
 =&\bigwedge_{\sigma,\tau\in\Sigma}\bigg([\sigma\ne\tau]\Rightarrow
 \bigvee_{\psi\in\mu(W)}[\sigma(\psi)\ne\tau(\psi)]\bigg)
 \\
 =&\mathbb{1},
 \end{align*} 
 where $c\Rightarrow d\!:=c^\perp\vee d$ for all $c,d\in \mathbb{B}$. By (1) $\mu(\widetilde{W})$ is wider then $\widetilde{\mu(W)}$. So that $[\![\mu(\widetilde{W})$  separates the points
of $\tilde{\Sigma}]\!]=\mathbb{1}$, as claimed.
 \end{proof}

 \medskip

 {\sc Proof of Theorem 3.8.}
 \begin{proof}
 There is no loss of generality in assuming that $\Lambda=\mathcal{R}{\downarrow}$ and  $X=\mathcal{X}{\downarrow}$,
 where $\mathcal{X}$ is a Banach space within $\mathbb{V}^{\mathbb{(\mathbb{}B)}}$. By
 Lemma 5.5, $[\![\widetilde{W}\subset C_{uo}(\tilde{S},\mathcal{X})]\!]=\mathbb{1}$ and
 $[\![\mu({\widetilde{W}})$ separates the points of $\tilde{S}]\!]=\mathbb{1}$.
 Take $f\in C_{uo}(\tilde{S},X)$ and note that $\tilde{f}\in C_{uo}(\tilde{S},\mathcal{X})$
 by Lemma 5.3. Using transfer  enables us to apply Theorem 1.1
 within $\mathbb{V}^{(\mathbb{B})}$, so that for ally
 $0<\varepsilon\in\mathbb{R}^{\scriptscriptstyle\wedge}$ the following are equivalent:

 $(1)$~There exists $\tilde{g}\in \widetilde{W}$ with
 $\|\tilde{f}-\tilde{g}\|<\varepsilon$.

 $(2)$~For every $s\in \tilde{S}$ there exists $\tilde{g}_s\in \widetilde{W}$ such that
 $\|\tilde{f}(s)-\tilde{g}_s(s)\|_{\mathcal{X}}<\varepsilon$.

 According to Lemma 5.3 there exist $g\in C_{uo}(S,X)$ such that
 $[\![\tilde{g}(u^{\scriptscriptstyle\wedge})=g(u)]\!]= \mathbb{1}$ for all $u\in S$.
 Clearly, $[\![\tilde{f}-\tilde{g}=(f-g)\,\tilde{{}}\,]\!]=\mathbb{1}$ and,
 $[\![\[f-g\]=\|\tilde{f}-\tilde{g}\|\,]\!]=\mathbb{1}$, so  $(1)$ is
 equivalent to $\inf\{\[f-g\]:\ g\in W\}=0$ in view of Lemma 4.4. Indeed, if
 $A\!:=\{\[f-g\]:\ g\in W\}$ and $B=\{\|\tilde{f}-\tilde{g}\,\|:\ \tilde{g}\in\widetilde{W}\}$,
 then for arbitrary $a\in\mathcal{R}{\downarrow}$ we have
 \begin{equation*}
 [\![a\in A{\uparrow}]\!]=\bigvee\{[\![a=u]\!]:\ u\in A\}
 \end{equation*}
 \begin{equation*}
 =\bigvee\{[\![a=\[f-g\]\,]\!]:\ g\in W\}
 =\bigvee\{[\![a=\|\tilde{f}-\tilde{g}\|\,]\!]:\ g\in W\}
 \end{equation*}
 \begin{equation*}
  =\bigvee\{[\![a=\|\tilde{f}-\tilde{g}\,\|:\
 \tilde{g}\in\widetilde{W}\}]\!]=[\![a\in B]\!],
 \end{equation*}
 hence
 $A{\uparrow}=B$. Taking into account Lemma 5.5, we can show in a similar
 way that $[\![\tilde{f}(s)-\tilde{g}(s)=(f-g)\,\tilde{{}}(s)\,]\!]=\mathbb{1}$ and,
 $[\![\[f(s)-g(s)\]=\|\tilde{f}(s)-\tilde{g}(s)\|\,]\!]=\mathbb{1}$, so the assertion $(2)$
 is equivalent to $\inf\{\[f(s)-g(s)\]:\ g\in W\}=0$ for all $s\in\Sigma$.~\end{proof}
 \medskip

 {\sc Proof of Theorem 3.12.}

 \begin{proof}
 Using the above notation, we see that $W$ and $f$ satisfy the conditions
 of Theorem 1.2 within $\mathbb{V}^{(\mathbb{B})}$. Define $A,B(\sigma)\in\mathbb{V}$
 and $\tilde{A},\tilde{B}(\sigma)\in\mathbb{V}^{(B)}$ by
 \begin{equation*}
 A\!:=\{\[f-g\]:\ g\in W\},\quad B(\sigma)\!:=\{\[f(\sigma)-x\]_{X}:\ x\in W(\sigma)\};
 \end{equation*}
 \begin{equation*}
 \tilde{A}\!:=\{\|\tilde{f}-g\|:\ g\in \tilde{W}\},\quad
 \tilde{B}(\sigma)\!:=\{\|\tilde{f}(\sigma)-x\|_\mathcal{X}:\ x\in\tilde{W}(\sigma)\}.
 \end{equation*}
 It can be shown as above that $[\![\tilde{A}=A{\uparrow}]\!]=\mathbb{1}$
 and $[\![\tilde{B}(\sigma)=B(\sigma){\uparrow}]\!]=\mathbb{1}$. By the
 Transfer Principle there exists $\sigma\in\tilde{\Sigma}$ such that
 $[\![\tilde{A}=\tilde{B}(\sigma)]\!]=\mathbb{1}$. It remains to apply Lemma 4.4 to obtain
 the required equality $A=B(\sigma)$.
 \end{proof}

 \Remark{~5.6}
 Theorems 3.8 and 3.12 contain Theorems 1.1 and 1.2  as particular cases. At the same time,
as seen from the above, Theorems 3.8 and 3.12 are the Boolean-valued interpretations of Theorems 1.1 and 1.2.
 Of course, the proofs avoiding Boolean valued analysis are also possible.

Using the same technique, we can formulated negative results on order approximation by
the cyclic hulls of sublattices and Grothendieck subspaces; cp. \cite[Ch.~3]{KKTop}.




 \Theorem{~5.7}~Let $E$ be a Dedekind complete vector lattice, while
 $S$ is a compact space, $\Lambda=\Orth(E)$, and $X=C_{uo}(S,E)$. Then

 $(1)$~a sublattice $L$ of $X$ is such that $\mix(L)$ does not coincide with $X$ if and only 
if there are $\Lambda$-linear lattice homomorphisms $S,T: X \to E$ satisfying $L\subset \ker(S-T)$;%

 $(2)$~A Grothendieck subspace $L$ of $X$ is such that
 $\mix(L)$ does not coincide with $X$ if and only 
if there are $\Lambda$-linear lattice homomorphisms $S,T: X \to E$  satisfying $L\subset\ker(S+T)$.
\Endproc

\section{Funding}
The research was supported by~the Ministry of Science and Higher Education of the Russian Federation (
Agreement 075--02--2022--896) and carried out in the framework of the State Task 
to the Sobolev Institute of Mathematics (Project 0314--2019--0005).

\medskip

\noindent
{\it Anatoly G.~Kusraev}\\
{\leftskip\parindent\small
\noindent
Southern Mathematical Institute\\
 22 Markus Street\\
Vladikavkaz, 362027, RUSSIA\\
E-mail: kusraev@smath.ru
\par}

\medskip
\noindent
{\it Sem\"en S.~Kutateladze}\\
{\leftskip\parindent\small
\noindent
Sobolev Institute of Mathematics\\
4 Koptyug Avenue\\
Novosibirsk, 630090, RUSSIA\\
E-mail: sskut@math.nsc.ru
\par}
\medskip

\end{document}